\documentclass[preprint]{imsart}
\RequirePackage[OT1]{fontenc}
\RequirePackage{amsthm,amsmath}
\RequirePackage[numbers]{natbib}
\RequirePackage[colorlinks,citecolor=blue,urlcolor=blue]{hyperref}

\issue{}
\firstpage{}
\lastpage{}

\startlocaldefs
\usepackage[tikz]{bclogo}
\usepackage{algorithm2e}
\usepackage[normalem]{ulem}
\usepackage{amssymb, amsmath, natbib, amsthm, graphicx, euscript, mathrsfs, enumerate }

\newcommand{\R}{\mathbb{R}}
\def\I{I\!\!I}
\renewcommand{\i}{\mathrm{i}}

\newcommand{\eqp}{\stackrel{p}{\longrightarrow}}
\newcommand{\eql}{\stackrel{d}{\longrightarrow}}

\renewcommand{\P}{{\mathbb P}}

\newcommand{\N}{{\mathbb N}}
\newcommand{\NN}{{\mathscr N}}

\newcommand{\E}{{\mathbb E}}

\newcommand{\eps}{\varepsilon}
\newcommand{\Var}{\operatorname{Var}}

\newcommand{\M}{\mathcal{M}}

\theoremstyle{plain}
\newtheorem{thm}{Theorem}[section]

\newtheorem{prop}[thm]{Proposition}
\theoremstyle{definition}

\newcommand{\wS}{\widetilde{S}}

\begin{document}
\begin{frontmatter}
\title{Limit theorems for sums of random variables with  mixture distribution}% \thanksref{T1}}
\runtitle{Limit theorems for  mixtures}

\begin{aug}

\author{\fnms{Vladimir} \snm{Panov}\ead[label=e1]{vpanov@hse.ru}}
\address{International Laboratory of Stochastic Analysis and its Applications\\National Research University Higher School of Economics \\ Shabolovka, 26, Moscow,  119049 Russia.\\
\printead{e1}}

%\thankstext{T1}{
%The study has been funded by the Russian Academic Excellence Project ``5-100''.}
\runauthor{V.Panov}

\end{aug}

\begin{abstract}
In this paper, we study the fluctuations of sums of random variables with distribution defined as a mixture of light-tail and truncated heavy-tail distributions. We focus on the case when both the mixing coefficient and the truncation level depend on the number of summands. The aim of this research is to characterize the limiting distributions of the sums due to various  relations between these parameters.
%these parameters which guarantee the fulfillment of the law of large numbers, central limit theorem, and convergence to other non-Gaussian distributions. 
\end{abstract}



\end{frontmatter}

\section{Introduction}
Theory of limit distributions for the sums of random variables is well-described in brilliant books by Ibragimov and Linnik \cite{il},   Meerschaert and Scheffler \cite{ms}, Petrov \cite{petrov}. Usually, the most interest is drawn to 2 classical models: a model of i.i.d. random variables and triangular arrays. For the first model,  it is common to find non-degenerate laws, which can appear as a limit of   the sums \((\xi_{1}+...+\xi_{n}) b^{-1}_{n} +a_{n}\) with i.i.d. \(\xi_{1},..., \xi_{n}\), and some deterministic sequences \(a_{n}, b_{n}\). It is well-known that the set of limiting distributions in this case coincides with the class of stable distributions.

In the second model, one considers an  infinitesimal triangular array - a collection of real random variables \(\{ Z_{nk}, k=1..k_{n}\},\) \(k_{n} \to \infty\) as \(n \to \infty,\)  such that \(Z_{n1}, ..., Z_{nk_{n}}\) are independent for each \(n\) and satisfy the condition of infinite smallness
\begin{eqnarray}
\label{inf}
\sup_{k=1..k_{n}}
\P \left\{
	\left|
		Z_{nk}
	\right| 
	> \delta
\right\} \to 0, \qquad n \to \infty
\end{eqnarray}
for any \(\delta>0.\) For this model, it is known  that only the infinitely divisible distributions can appear as the non-degenerate limit of sums \(Z_{n1}+...+Z_{n k_{n}} -a_{n}\) with deterministic \(a_{n}\), and moreover, for any infinitely-divisible distribution there exists a triangular array \((Z_{nk})\) such that the sum \(Z_{n1}+...+Z_{n k_{n}}\) converges to this distribution. 

Nevertheless, the analysis of the limiting distribution in particular models can be rather tricky. For instance,  Ben Arous, Bogachev and Mochanov \cite{benarous} analyzed the asymptotic behaviour of  the sums \(\sum_{i=1}^{n}e^{t \xi_{i}}\), where \(\xi_{1},...,\xi_{n}\) is an i.i.d. sequence of r.v. with regularly varying log-tail functions, and \(n\) and \(t\) simultaneously tend to infinity, provided that  
the speeds of growth of \(n\) and \(t\) are coordinated via the parameter 
\(\lambda=\liminf_{n,t \to \infty} \log(n) / \log(\E[e^{t\xi}])\).
It turns out, that there exist two critical values of this parameter, \(\lambda_{1}<\lambda_{2}\),  below which the law of large numbers and the central limit theorem (respectively) break down.  This result generalizes some previous findings related to the random energy model, which corresponds to the case when \(\xi_{i}\) are standard normal, see Bovier, Kurkova and L{\"o}we \cite{BKL}.

In the current research, we derive similar results for  completely another model,  defined as a mixture of two distributions:  the first distribution has light tails and the second  is constructed by truncation of the distribution with heavy tails. More precisely,   let \(F_{1}(\cdot)\) be a distribution function corresponding to a probability distribution on \(\R_{+}\) with  the upper tail in exponential form, that is,
\begin{eqnarray*}
 F_{1}(x) = 1 - e^{-\lambda x} \left(
	1 +o(1)
\right), \qquad x \to +\infty, \end{eqnarray*}
for some \(\lambda>0\). %Assume also that \(\P_{1}\) has the central absolute moment  of order larger than 2:
%\begin{eqnarray}
%\label{A}
%\exists \delta>0: \qquad \int_{\R_{+}}
%\left|
%	x - \mu
%\right|^{2+\delta} \P_{1} (dx)<\infty \qquad \mbox{with} \quad
%\mu=\int_{\R_{+}}
 %y P_{1}(dy).
%\end{eqnarray}
Let \(F_{2}(\cdot)\)  be a distribution function corresponding to a heavy-tailed distribution with support on \([1,\infty)\), 
\begin{equation*}%\label{B}
F_{2}(x)
 = 1 - x^{-\alpha} \left(
1 + 
o(1)
\right), \qquad x \to +\infty,
\end{equation*}
with \(\alpha \in (0,2). \)  
By \(F_{2}^{M}(x) \) denote the truncated  distribution function  \(F_{2}\) at  level \(M\):
\begin{eqnarray*}
%\label{F2M}
F_{2}^{M} (x)  =
 \begin{cases}
		F_{2}(x)/F_{2}(M), &\text{if $x\leq M$;}\\
		 0, &\text{if $x> M$}.
	\end{cases}
\end{eqnarray*}
 Next, consider the mixture of these distributions, that is, the distribution with distribution function
\begin{eqnarray}
\label{F}
F(x) = (1 - \eps) F_{1}(x) + \eps F_{2}^{M}(x),
\end{eqnarray}
where \(\eps \in (0,1)\)  is a mixing parameter, which is assumed to be small.  The motivation of considering such mixture goes to the idea to model some ``frequent events'' by light-tail distributions, and ``rare events'' by truncated heavy-tailed. For instance, this idea is quite natural for modelling the claim amounts in insurance, see e.g, Rolski et al. \cite{RSST}, or Embrechts, Kl\"uppelberg and Mikosch \cite{EKM}.
Another example comes from  population dynamics by analyzing the migration of species. In some models, it is assumed that   for most species in the population, the distribution of migration is light-tailed, whereas for some small amount of species the distribution is heavy-tailed, see Whitmeyer and Yang~\cite{WY}.

In this research we consider the case when \(\eps\) and \(M\) depend on \(n\), and moreover,  \(\eps \to 0\) and \(M\to \infty\) as \(n\) grows. We focus on studying the fluctuations of sums of random variables drawn from the mixture model~\eqref{F}, and aim to characterize the limit laws  depending on the relation between \(\eps\) and \(M.\) Several problems of this type are considered in the paper by Grabchak and Molchanov \cite{gm}, where the parameter \(M\) as well as  both distributions in~\eqref{F}, are fixed. As it is shown in 
Grabchak and Molchanov \cite{gm}, the complete asymptotic analysis can be done by taking into account that the distributions  are in the domain of attractions of some stable random variables. In this case, it is clear that the limit law for the mixture can be determined by the relation between the normalizing sequences.  Nevertheless, this methodology cannot be applied to our set-up, since the parameters \(\eps\) and \(M\) simultaneously vary.

%The study of the asymptotic properties of this mixtures depending on the relations between \(\eps\) and \(M,\)
% is quite important for statistical methods. For instance, fluctuations of sums play key role for statistical approaches used in insurance and finance, see Embrechts, Kl\"uppelberg and Mikosch \cite{EKM}. 

The paper is organized as follows. In the next section we formulate our main results. It turns out (and  is not surprising) that the cases \(\alpha\in(0,1)\)  and \(\alpha \in [1,2)\) are essentially different, see Subsections~\ref{sec21} and~\ref{sec22} respectively. The proofs are collected in Section~\ref{proofs}.

\section{Limit theorems }
Assume that for any \(n \in \N\) we are given by  \(k_{n}=n\) independent  random variables \(Z_{n1},..., Z_{nk_{n}}\) with  mixing distribution  \eqref{F}.  In other words,\begin{eqnarray}
\label{mix}
Z_{nk} =  \left( 1 - B_{nk} \right) X_{nk} + B_{nk} Y_{nk}, \quad k=1..k_{n},
\end{eqnarray}
where \(X_{n1}, ..., X_{nk_{n}} \sim F_{1}, \;  Y_{n1}, ..., Y_{n k_{n}} \sim F_{2}^{M_{n}}, \; \) 
\(B_{n 1}, ..., B_{n k_{n}}\) are  Bernoulli random variables with probability of success equal to \(\eps_{n},\)  and all \(X_{nk}, Y_{nk}, B_{nk},\) \(k=1..k_{n}\) are jointly independent for any \(n\). 

In what follows, we take  \(M=n^{\gamma_{1}}, \; \eps = n^{-\gamma_{2}}\) with positive \(\gamma_{1},\gamma_{2},\) and aim to characterize the asymptotic behaviour of the sum \(S_{n}:=\sum_{k=1}^{k_{n}} Z_{nk}\) due to the relation between \(\gamma_{1}\) and \(\gamma_{2}.\) 

 \subsection{Case $\alpha \in (0,1)$}\label{sec21}
We start with the most interesting case, \(\alpha \in (0,1)\).
%We start with the law of large numbers.
\begin{thm}\label{thm3}

\begin{enumerate}[(i)]
\item
Let \(\gamma_{1}, \gamma_{2}\) be such that 
\[ \gamma_{2}>(2-\alpha) \gamma_{1} \qquad \mbox{or} \qquad  \gamma_{2}<\min \Bigl\{ 
(2-\alpha) \gamma_{1}, 1-\alpha \gamma_{1} 
\Bigr\}.\]
Then the central limit theorem holds, in the sense that 
\begin{eqnarray*}
	\frac{S_{n} - n \E\left[ Z_{n1}\right]
	}
	{
		\sqrt{n \Var(Z_{n1})}
	}
	\eql \NN(0,1), \qquad n \to \infty.
\end{eqnarray*}
\item  Let \(\gamma_{1}, \gamma_{2}\) be such that 
\begin{eqnarray*}
\gamma_{1}>1/2 \qquad \mbox{and}
\qquad
\gamma_{2} \in \left(
	1-\alpha/2, \; (2-\alpha) \gamma_{1}
\right).
\end{eqnarray*}
Then 
\begin{eqnarray*}
\frac{ 
S_{n}  - n \E [X_{n1}]
}
{
\sqrt{n \Var (X_{n1})}
}
	\eql \NN(0,1), \qquad n \to \infty.
\end{eqnarray*}

\item Finally, let \(\gamma_{1}, \gamma_{2}\) be such that 
\begin{eqnarray*}
\gamma_{1}>1/2 \qquad \mbox{and}
\qquad
 \gamma_{2} \in \left(
	\max\left(
		1-\alpha \gamma_{1}, 0
	\right),
	1 - \alpha/2
\right).
\end{eqnarray*}
Then for any constant \(c>0\)
\begin{eqnarray}
\label{convtostab}
	\frac{S_{n}-  \kappa_{n}}{c n^{(1-\gamma_{2})/\alpha}}  \eql 
	F_{\alpha, c}, \qquad n \to \infty,
\end{eqnarray}
where \(F_{\alpha, c}\) is an \(\alpha\)-stable distribution on \(\R_{+}\), that is, an infinitely divisible distribution with the L{\'e}vy  density  \(s(x) = c x^{-1-\alpha} \I\{ x >0\}\),  and without continuous part,
and 
\begin{eqnarray}\label{kappan}
	\kappa_{n} = \begin{cases}
		n\E[X_{n1}], &\text{if $\gamma_{2}\in (1-\alpha,1-\alpha/2)$,}
		\\
		n\E[X_{n1}] +
		\frac{c \alpha}{1-\alpha} n^{(1-\gamma_{2})/\alpha}, &\text{if $\gamma_{2}=1-\alpha$,}\\
		\frac{c \alpha}{1-\alpha} n^{(1-\gamma_{2})/\alpha}, &\text{if $\gamma_{2}\in  (0, 1-\alpha)$.}
	\end{cases}
\end{eqnarray}

%\[s(x) = \frac{c_{1}}{x^{1+\alpha}} \I\{x>0\} + \frac{c_{2}}{|x|^{1+\alpha}} \I\{x<0\} \]
%with \(c_{1} \geq 0, c_{2} \geq 0, c_{1}+c_{2}>0.\)
\end{enumerate}
\end{thm}

Note that the normalizing term % in (i) can be changed to \(\sqrt{n \Var(X_{n1})}\), but normalizing term 
in (ii) cannot be changed to \(\sqrt{n \Var(Z_{n1})}\), and therefore (ii) essentially differs from the central limit theorem. In fact, 
\begin{eqnarray}\label{var}
\Var(Z_{n1})=\left(
\Var(X_{n1}) + \frac{2 \alpha}{2- \alpha} n^{(2-\alpha) \gamma_{1} - \gamma_{2}}
\right) (1+ o(1)), \quad n \to \infty\end{eqnarray}
and hence \(\Var(Z_{n1})\asymp\Var(X_{n1})\) if and only if \(\gamma_{2}>(2-\alpha) \gamma_{1}.\)

\begin{thm}\label{thm1}
\begin{enumerate}[(i)]
\item 
Let \(\gamma_{1}, \gamma_{2}\) be such that 
%\begin{enumerate}[1.]
\[
 \gamma_{2}>(1-\alpha) \gamma_{1}\qquad \mbox{or}\qquad\gamma_{2}<\min \Bigl\{ 
(1-\alpha) \gamma_{1}, 1-\alpha \gamma_{1} 
\Bigr\}.\]
Then the law of large numbers holds, in the sense that 
\begin{eqnarray*}
	\frac{S_{n}}{n \E\left[Z_{n1}\right]} \eqp 1, \qquad n\to \infty.
\end{eqnarray*}
\item Let \(\gamma_{1}, \gamma_{2}\) be such that 
\begin{eqnarray*}
\gamma_{2}>1-\alpha \qquad \mbox{and}
\qquad
\gamma_{2} \in \left(
	1-\alpha, (1-\alpha )\gamma_{1}
\right).
\end{eqnarray*}
Then the analogue of the law of large numbers with normalization \(n \E[X_{n1}]\) holds, i.e.,
\begin{eqnarray*}
	\frac{S_{n}}{n \E\left[X_{n1}\right]} \eqp 1, \qquad n\to \infty.
\end{eqnarray*}
%\item Let \(\gamma_{2} \in \left(
%	\max\left(
%		1-\alpha \gamma_{1}, 0
%	\right),
%	1 - \alpha
%\right).
%\) Then 
%
%\begin{eqnarray*}
%	S_{n} \eqp 1, \qquad n\to \infty.
%\end{eqnarray*}
\end{enumerate}
\end{thm}
Analogously to \eqref{var}, we note that 
\begin{eqnarray}
\label{EZZ}
\E\left[
Z_{n1}
\right]
=
\left(
\E\left[
X_{n1}
\right] +\frac{\alpha}{1-\alpha} n^{(1-\alpha)\gamma_{1} - \gamma_{2}}
\right)(1+o(1)), \quad n \to \infty,
\end{eqnarray}
and therefore \(\E\left[
Z_{n1}
\right]
\asymp
\E\left[
X_{n1}
\right] \) if and only if \(\gamma_{2}>(1-\alpha) \gamma_{1}.\)

\begin{figure}
\begin{center}
\includegraphics[width=0.6\linewidth ]{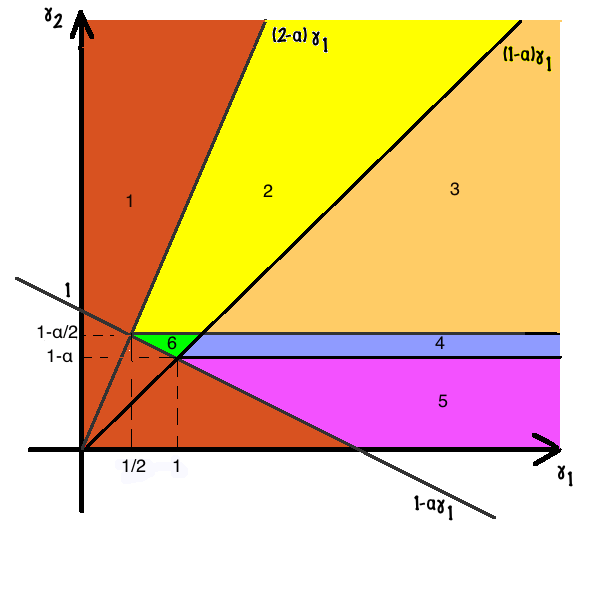}\caption{\label{fig1}Illustration of the  limit behavior of the sums \(\sum_{k=1}^{n} X_{nk}\) depending on \(\gamma_{1}\) and \(\gamma_{2}\) for the case \(\alpha \in (0,1)\). \newline zone 1(red): both  the central limit theorem and the law of large numbers  hold (see Theorem~\ref{thm3}(i) and Theorem~\ref{thm1}(i)); \newline
zone 2 (yellow): convergence to the standard normal distribution under normalization \(n \E[X_{n1}]\) and the law of large numbers (see Theorem~\ref{thm3}(ii) and Theorem~\ref{thm1}(i));\newline
zone 3 (orange): convergence to the standard normal distribution under normalization \(n \E[X_{n1}]\) and the analogue of the LLN with the same normalization (see Theorem~\ref{thm3}(ii) and Theorem~\ref{thm1}(ii));\newline 
zone 4 (blue):  convergence to stable distribution and the analogue of the LLN with normalization \(n \E[X_{n1}]\) (see Theorem~\ref{thm3}(iii) and Theorem~\ref{thm1}(ii)); \newline 
zone 5 (purple): convergence to stable distribution (see Theorem~\ref{thm3}(iii));\newline
zone 6 (green): convergence to stable distribution and LLN  (see Theorem~\ref{thm3}(iii) and Theorem~\ref{thm1}(i)).}
\end{center}
\end{figure}

Figure~\ref{fig1} illustrates the division of the area \(\left(
	\gamma_{1}, \gamma_{2}
\right) \in \R_{+} \times \R_{+}\) into subareas with different asymptotic properties of the sums \(\sum_{k=1}^{n} X_{nk}.\)
\subsection{Case $\alpha \in [1,2)$}
\label{sec22}
In this case, the limit laws are more simple. We formulate the results in the next theorem. 
\begin{thm}\label{other}
Let \(\alpha \in [1,2).\) Then 
\begin{enumerate}[(i)]
\item for any positive \(\gamma_{1}, \gamma_{2}\) 
the law of large numbers  holds, i.e.,
\begin{eqnarray*}
	\frac{S_{n}}{n \E\left[Z_{n1}\right]} \eqp 1, \qquad n\to \infty;
\end{eqnarray*}

\item if \(\gamma_{1}, \gamma_{2}\) are such that 
\begin{eqnarray}
\label{cond} \gamma_{2}>(2-\alpha) \gamma_{1} \qquad \mbox{or} \qquad  \gamma_{2}<\min \Bigl\{ 
(2-\alpha) \gamma_{1}, 1-\alpha \gamma_{1} 
\Bigr\},
\end{eqnarray}
then the central limit theorem holds, i.e.,
\begin{eqnarray*}
	\frac{S_{n} - n \E\left[ Z_{n1}\right]
	}
	{
		\sqrt{n \Var(Z_{n1})}
	}
	\eql \NN(0,1), \qquad n \to \infty;
\end{eqnarray*}
\item otherwise, if \eqref{cond} is not fulfilled, then for any constant \(c>0,\)
\begin{eqnarray*}
	\frac{S_{n} - n \E\left[ Z_{n1}\right]
	}
	{
		c n^{(1-\gamma_{2})/\alpha}
	}
	\eql F_{\alpha,c}, \qquad n \to \infty,
\end{eqnarray*}
see Theorem~\ref{thm3}(iii) for notations.
\end{enumerate}
\end{thm}

\section{Proofs}\label{proofs}
We first prove the statements related to the law of large numbers  (Theorem~\ref{thm1} (i) and Theorem~\ref{other} (i)), then the central limit theorems (Theorem~\ref{thm3}(i) and Theorem~\ref{other}(ii)), and afterwards we show the convergence to stable non-Gaussian distributions (Theorem~\ref{thm3} (ii), (iii) and Theorem~\ref{other}(iii)).

%Below we consider the case \(\alpha \in (0,1) \) in details; the proof for all other cases follows the same lines.
\textbf{Proof of Theorem~\ref{thm1} (i) and Theorem~\ref{other} (i).} 

\textbf{1.} Denote 
\begin{eqnarray*}
\wS_{n}=\frac{S_{n}}{\E[S_{n}]} -1 =
 \frac{
 	\sum_{k=1}^{n}
	\left(
		Z_{n k} - \E \left[
			Z_{nk}
		\right]
	\right)
}{
	n\E \left[
		Z_{n1}
	\right]
}.
\end{eqnarray*}
Our aim is to show that there exists a constant \(r>1\) such that \(\E
	| 
		\wS_{n}
	|^{r} 
 \to 0\) as \(n \to \infty\). This will imply that \(\wS_{n} \eqp 0,\) and therefore the result will follow.

Applying the Bahr-Esseen inequality for \(r \in(1,2),\)  see \cite{bahresseen}, we get that 
\begin{eqnarray}
\nonumber
\E\left[
	| 
		\wS_{n}
	|^{r} 
\right]
&\leq&
C_r
\frac{
 	\sum_{k=1}^{n}
	\E
	\left[
		|
			Z_{n k} - \E \left[
					Z_{nk}
				\right]
		|^{r}
	\right]
}{
	\left( n \E \left[
		Z_{n1}
	\right]
	\right)^{r}
}\\
\label{BE}
&=&
C_r n^{1-r} \frac{\M_{n}(r)}{\left( \E \left[
		Z_{n1}
	\right]
	\right)^{r}},
\end{eqnarray}
where \(C_r\) is  some constant depending on \(r\), and \(\M_{n}(r)\) is the \(r\)-th absolute central moment of \(Z_{n1}\). The further analysis is consists in establishing the asymptotical behavior of the numerator and denominator of the fraction in \eqref{BE}.

\textbf{2.} Note that 
 for any \(n=1,2,..\)
\begin{eqnarray*}
\E [Z_{n1}]&=&%\E [ 1 - B_{n1} ] \cdot \E[ X_{n1} ]+ \E [B_{n1}]\cdot\E[ Y_{n1}]=
(1-\eps_{n}) \cdot \mu_{1}(1) + \eps_{n} \cdot \mu_{2}(1),
\end{eqnarray*}
where \(\mu_{1}(s) = \E[(X_{n1})^{s}] < \infty\) for any \(s>0\), and 
 \begin{eqnarray*}
\mu_{2}(s)  =
\E[(Y_{n1})^{s}] 
=
\begin{cases}
 \frac{\alpha}{|s-\alpha|} M_{n}^{\max(s-\alpha,0)} \left(
	1 + o(1)
\right),  &\text{if $s \ne \alpha$;}\\
		\alpha \log(M_{n})\left(
	1 + o(1)
\right), &\text{if $s=\alpha$.}
\end{cases}%, \qquad n \to \infty.
\end{eqnarray*}
 Therefore, 
\begin{eqnarray}\label{EZ}
	\E [Z_{n1}] &=&
	\begin{cases}
\bigl(
		\mu_{1}(1) +
		\frac{\alpha}{1-\alpha} \eps_{n}M_{n}^{1-\alpha} 
	\bigr) \left(
		1 + o(1)
	\right),  &\text{if $\alpha<1$;}\\
	\bigl(
		\mu_{1}(1) +
		\alpha \eps_{n}\log(M_{n})
	\bigr)  \left(
	1 + o(1)
\right), &\text{if $\alpha=1$;}\\
		\mu_{1}(1)\left(
	1 + o(1)
\right), &\text{if $\alpha>1$.}
\end{cases}
\end{eqnarray}

\textbf{3.} It holds for any \(r>1\)
\begin{eqnarray*}
\M_{n}(r)  &=& 
\E\Biggl[
\biggl|
	 \left( 1 - B_{n1} \right) X_{n1} + B_{n1} Y_{n1}
	 - 
	 \E \Bigl[
		 \left( 1 - B_{n1} \right) X_{n1} + B_{n1} Y_{n1}	 	
	 \Bigr]
\biggr|^{r}
\Biggr]\\
&\leq & 4^{r-1} \Biggl\{
\E\left[
	(1-B_{n1})^{r}
\right]
\cdot
\E\left[
	X_{n1}^{r}
\right]
+
\E\left[
	B_{n1}^{r}
\right]
\cdot
\E\left[
	Y_{n1}^{r}
\right] \Biggr.
\\
&&
\Biggl.
\hspace{4cm}+\left(
	1-\eps_{n}
\right)^{r}
\cdot
\left(
	\E X_{n1}
\right)^{r}
+ \eps_{n}^{r} 
\cdot
\left(
	\E Y_{n1}
\right)^{r}
\Biggr\}\\
&=&
4^{r-1}
\Biggl\{
	(1-\eps_{n}) \mu_{1}(r) +\eps_{n} \mu_{2}(r) + 
	(1-\eps_{n})^{r} \left(
		\mu_{1}(1) 
	\right)^{r}
	\\
&& \hspace{7cm}
	+\eps_{n}^{r}
	\left(
		\mu_{2}(1) 
	\right)^{r}
\Biggr\}.
\end{eqnarray*}
Denote \(D_{r}= 
	\mu_{1}(r) + (\mu_{1}(1))^{r}\) and consider two cases:
\begin{enumerate}[(a)]
\item if \((\alpha<1, r>1)\) or \((1<\alpha\leq r)\), then it holds
\begin{eqnarray*}
	\M_{n}(r)  &\leq &
	4^{r-1}
 \left( 
D_{r} +\frac{\alpha}{r-\alpha} \eps_{n} M_{n}^{r- \alpha}
\right) 
\cdot 
\left(
	1  + o(1)
\right),
\end{eqnarray*}
where  we use that \(\eps_{n} \mu_{2}(r) \gtrsim \eps_{n}^{r}
	\left(
		\mu_{2}(1) 
	\right)^{r}\) as \(n \to \infty;\) 
\item
otherwise, if $(1<r<\alpha)$ or \((\alpha=1, r>1)\), then we have 
\begin{eqnarray*}
	\M_{n}(r)  &\leq &
	4^{r-1}
D_{r}
\left(
	1  + o(1)
\right).\end{eqnarray*}
\end{enumerate}
\textbf{4.} To conclude the proof, we substitute the upper estimate for \(\M_{n}(r)\) and \eqref{EZ} into \eqref{BE}.
If \(\alpha<1,\) then 
\begin{eqnarray}
\label{ES}
\E\left[
	| 
		\wS_{n}
	|^{r} 
\right]
&\leq&
\widetilde{C}_{r}\cdot 
n^{1-r} 
\frac{
	\left( 	
		c_{r}^{(1)} + \eps_{n} M_{n}^{r- \alpha} 
	\right) 
	\cdot
	\left(
		1  + o(1)
	\right)
}{
	\left( 	
		c_{r}^{(2)}+ \eps_{n} M_{n}^{1- \alpha} 
	\right)^{r} 
	\cdot
	\left(
		1  + o(1)
	\right)
}
\end{eqnarray}
with some  constants \(c_{r}^{(1)},c_{r}^{(2)}\) (depending on \(r\)) and a bounded function \(\widetilde{C}_{r}.\) The asymptotic behaviour of the last fraction differs between the following two cases:
\begin{itemize}
\item \(\eps_{n} M_{n}^{1- \alpha}  \to 0,\) that is, \(\gamma_{2}>(1-\alpha) \gamma_{1}\). Then there exists \(r>1\) such that \(\eps_{n} M_{n}^{r- \alpha}  \to 0\)  - in fact, one can take \(r=(\gamma_{2} - (1-\alpha)\gamma_{1})\gamma_{1}^{-1}+1>1\). Under this choice of \(r,\) the r.h.s. of \eqref{ES} tends to \(0,\) and therefore the law of large numbers holds for any \((\gamma_{1}, \gamma_{2})\) s.t. \(\gamma_{2}>(1-\alpha) \gamma_{1}\).
\item \(\eps_{n} M_{n}^{1- \alpha}  \to \infty,\) that is, \(\gamma_{2}<(1-\alpha) \gamma_{1}.\) Then the right-hand side of \eqref{ES} tends to 0 if and only if \(M_{n}^{\alpha} /( n \eps_{n}) \to 0\), that is, \(\gamma_{2}<1-\alpha \gamma_{1}.\) This case corresponds to the area \(\gamma_{2}<\min \Bigl\{ 
(1-\alpha) \gamma_{1}, 1-\alpha \gamma_{1} 
\Bigr\}.\)
 \end{itemize}
 
In other cases, \(\alpha>1\) and \(\alpha=1,\) we can choose \(r \in (1, \alpha)\) and get that \(\E\left[
	| 
		\wS_{n}
	|^{r} 
\right] \lesssim n^{1-r},\) and therefore the law of large numbers holds with any positive \(\gamma_{1},\gamma_{2}\).

\textbf{Proof of Theorem~\ref{thm3}(i) and Theorem~\ref{other}(ii)}
To prove these theorems, we check that the Lyapounov condition holds (see (27.16) from \cite{Bill}): there exists \(\delta>0\) such that
\begin{eqnarray*}
\Omega_{n}:=
\frac{\M_{n}(2+\delta)}{
	n^{\delta/2} \left(
		\Var(Z_{n1})
	\right)^{1+\delta/2}
} \to 0, \qquad \mbox{as} \quad 
n \to \infty.
\end{eqnarray*}
The variance of \(Z_{n1}\) has the following asymptotical behaviour:
\begin{eqnarray*}
\Var(Z_{n1}) &=& 
\left(
	1-\eps_{n}
\right)
 \mu_{1}(2) 
 + \eps_{n} \mu_{2}(2) - \left(
\E [Z_{n1}]
\right)^{2}
%\\
%&=& 
%\left(
%	1-\eps_{n}
%\right)
% \mu_{1}(2) 
% +
% \frac{\alpha}{2-\alpha} \eps_{n} M_{n}^{2-\alpha} (1 +o(1))
% \\
% && 
% \hspace{2cm}
%- \left\{
%	\left(
%		1 - \eps_{n}
%	\right)	 \mu_{1}(1) + 
%	 \frac{\alpha}{1-\alpha} \eps_{n}M_{n}^{1-\alpha} (1 +o(1))
%\right\}^{2}
\\
&=&	\biggl(
		\left[ \mu_{1}(2) - (\mu_{1}(1))^{2}\right] +
		\frac{\alpha}{2-\alpha} \eps_{n}M_{n}^{2-\alpha} 
	\biggr) \cdot \left(
		1 + o(1)
	\right), \qquad n \to \infty,
\end{eqnarray*}
and the numerator of \(\Omega_{n}\) was already considered in the proof of Theorem~\ref{thm1} (i). Therefore, 
\begin{eqnarray*}
\Omega_{n} \leq
c_{1} \cdot
\frac{
	\left( 	
		c_{2} + \eps_{n} M_{n}^{2+\delta- \alpha} 
	\right) 
	\cdot
	\left(
		1  + o(1)
	\right)
}{
	n^{\delta/2}
	\left( 	
		c_{3} + \eps_{n} M_{n}^{2- \alpha} 
	\right)^{1+\delta/2} 
	\cdot
	\left(
		1  + o(1)
	\right)
}, \qquad n \to \infty,
\end{eqnarray*}
with some positive constants \(c_{1}, c_{2}, c_{3}.\)
The rest of the proof follows the same lines as Step~4 in the proof of Theorem~\ref{thm1}(i), see above.

\textbf{Proof of Theorem~\ref{thm3} (ii), (iii) and Theorem~\ref{other}(iii).} The proof is based on the following proposition, which is in fact a combination of  Theorem~1.7.3 from \cite{il}, Theorem~3.2.2 from \cite{ms}, and a number of theorems given in Chapter~IV from \cite{petrov}.
\begin{prop}
\label{con}
Consider an infinitesimal triangular array  \(\{ Z_{nk}, k=1..k_{n}\},\) such that \eqref{inf} is fulfilled.  In what follows, we denote the distribution of \(Z_{nk}\) by \(\mu_{nk},\) and use the notation 
\(S_{n}:=Z_{n1}+... +  Z_{nk_{n}}\). The following statements hold.
\begin{enumerate}
\item  If there exists a random variable \(Y\) and a  sequence of real numbers \(a_{n} \) such that 
\begin{eqnarray}
\label{ZZ}
S_{n} - a_{n} \eql Y, \qquad n \to \infty, 
\end{eqnarray}
then \(Y\) has an infinitely divisible distribution; moreover, for any infinitely distribution \(\P_{inf}\) there exists a triangular array \(\{ Z_{nk}, k=1..k_{n}\}\) such that \(S_{n}\eql \P_{inf}\).
\item There exists a deterministic sequence \(a_{n}\) such that sequence \(S_{n} - a_{n}\)  converges weakly to an infinitely divisible random variable \(Y\) with 
characteristic exponent
\begin{eqnarray*}
\psi(u) = \i u \mu - \frac{1}{2} u^{2} \sigma^{2} +
\int_{\R/\{0\}}\left(
	e^{\i u x} - 1 - \i u x \I\left\{
	|x| \leq 1
\right\}
\right) \nu(dx),
\end{eqnarray*}
where  \((\mu, \sigma, \nu)\) is a L{\'e}vy triplet,  if and only if the following conditions are fulfilled:
\begin{enumerate}
\item \( \sum_{k=1}^{k_{n}}\mu_{nk} (A)\to \nu\left(
	A
\right)\)
for any \(A=(-\infty,x)\) with \(x<0\) and any 
\(A=(x, +\infty)\) with \(x>0\) such that \(\nu(\partial A)=0\);
\item moreover,
\begin{multline}
\label{condB}
\lim_{\tau\to 0} \limsup_{n \to \infty}\sum_{k=1}^{k_{n}} \left\{
	\int_{|x|<\tau} x^{2} \mu_{nk} (dx) - \left(
		\int_{|x|<\tau} x \; \mu_{nk}(dx)
\right)^{2}
\right\}
\\
=\lim_{\tau\to 0} \liminf_{n \to \infty}\sum_{k=1}^{k_{n}} \left\{
	\int_{|x|<\tau} x^{2} \mu_{nk} (dx) - \left(
		\int_{|x|<\tau} x \; \mu_{nk}(dx)	
\right)^{2} \right\}= \sigma^{2}
\end{multline}
\end{enumerate}
If these conditions are satisfied, \(a_{n}\) may be chosen according to the formula 
\begin{eqnarray}\label{an}
	a_{n}= \sum_{k=1}^{k_{n}} \int_{|x|<1} x\; \mu_{nk}(dx) + o(1),
\end{eqnarray}
provided \(\nu\left\{
	x: \; |x| = 1
\right\}=0\). %In particular, if \(Y\) has a stable distribution then \(H\) is an arbitrary positive number.
\item There exists a deterministic sequence \(a_{n}\) such that sequence \(S_{n} - a_{n}\)  converges weakly to a standard normal random variable \(Y\) if and only if the following conditions are fulfilled: 
\begin{enumerate}
\item \(\sum_{k=1}^{k_{n}} \P \left\{
	\left|
		Z_{nk}
	\right| 
	>x
\right\} \to 0\) as  \(n \to \infty\) for any \(x>0\);% (note also that this condition guarantees \eqref{inf});
\item \(
 \lim_{n \to \infty}\sum_{k=1}^{k_{n}} \left\{
	\int_{|x|<\tau} x^{2} \mu_{nk} (dx) - \left(
		\int_{|x|<\tau} \mu_{nk} (dx)
\right)^{2} \right\}= 1
\)
for some \(\tau>0\).
\end{enumerate}
If these conditions are satisfied, \(a_{n}\) may be chosen according to \eqref{an}.
\end{enumerate}
\end{prop}
Returning to our setup, we denote \(F_{nk}(x) = \P\left\{Z_{nk} \leq \beta_{n} x \right\}\), and first note that for any \(x \in (\beta_{n}^{-1}, \beta_{n}^{-1} M_{n})\)
\begin{eqnarray*}
\sum_{k=1}^{k_{n}} \left(
1 - F_{nk}(x)
\right)&=& 
\sum_{k=1}^{k_{n}}
\left[
	1 - \left(
		1 - \eps_{n}
	\right) F_{1} (\beta_{n} x)
	-
	\eps_{n}\frac{F_{2} (\beta_{n} x)}{F_{2}(M_{n})}
\right]\\
&=& 
n (1 - \eps_{n}) \left(
	1 -F_{1} (\beta_{n} x)
\right)
+ n \eps_{n}  \left(
	1 -\frac{F_{2} (\beta_{n} x)}{F_{2}(M_{n})}
\right)\\
&=& n (1 - \eps_{n}) e^{-\lambda \beta_{n} x}
\left(
	1 + o(1)
\right)
+ n \eps_{n}  \left(
\beta_{n} x
\right)^{-\alpha}
\left(
	1 + o(1)
\right)\\
&& \hspace{5cm}-  n \eps_{n}  M_{n}^{-\alpha}
\left(
	1 + o(1)
\right).
\end{eqnarray*}
Note that basically only 3 situations are possible. 
\begin{enumerate}[1.] 
\item \underline{\(1-\alpha \gamma_{1}<\gamma_{2}<1.\)} In this case, under the choice \(\beta_{n}=c_{1} n^{(1-\gamma_{2})/\alpha}\) with any constant \(c_{1}>0\) we get 
\begin{eqnarray*}
\sum_{k=1}^{k_{n}} \left(
1 - F_{nk}(x)
\right) \to c_{1} x^{-\alpha}, \qquad \forall \; x \in \R_{+},
\end{eqnarray*}
because \(n e^{-\lambda \beta_{n}} \to 0\), \(n\eps_{n}\beta_{n}^{-\alpha} \to c_{1}\), and \(n \eps_{n} M_{n}^{\alpha } \to 0\). Moreover, the condition \eqref{inf} is fulfilled - in fact, for any \(\delta>0,\) it holds 
\begin{eqnarray*}
\sup_{k=1..k_{n}}
\P \left\{
	\left|
		Z_{nk}
	\right| 
	> \delta
\right\} 
=k_{n}^{-1}\sum_{k=1}^{k_{n}} \left(
1 - F_{nk}(\delta)
\right) \to 0.
\end{eqnarray*}
 Next, with any \(s\geq 1,\)
\begin{eqnarray*}
\int_{|x|<\tau} x^{s} \; \widetilde\P_{1} (dx) &=& \beta_{n}^{-s}\cdot  \E \Bigl[
	X_{n1}^{s} \Bigr]  \left( 1 + o(1) \right), \\
\int_{|x|<\tau} x^{s} \; \widetilde\P_{2} (dx)
&=&
	\begin{cases}
\beta_{n}^{-\min(\alpha,s)}   \tau^{\max(s-\alpha,0)} \frac{\alpha}{|s-\alpha|} \left( 1 + o(1) \right),
&\text{if $\alpha\ne s$;}\\
\alpha \beta_{n}^{-s}  \log(\beta_{n}) \left( 1 + o(1) \right),
&\text{if $\alpha = s$.}
\end{cases}
\end{eqnarray*}
where \(\widetilde{\P}_{1},\widetilde{\P}_{2}\) are the probability distributions of \(X_{n1} / \beta_{n}\) and \(Y_{n1} / \beta_{n}\) resp.
Therefore, if \(\alpha<1\), the condition \eqref{condB} reads as 
\begin{multline*}
G_{n}:=
\sum_{k=1}^{k_{n}} \left\{
	\int_{|x|<\tau} x^{2} \mu_{nk} (dx) - \left(
		\int_{|x|<\tau} x \; \mu_{nk}(dx)	
\right)^{2} \right\}\\
= 
n \Biggl\{
	\beta_{n}^{-2} \E\left[
		X_{n1}^{2}
	\right]
	+ 
	C_{1}\eps_{n} \beta_{n}^{-\alpha}\tau^{2-\alpha} 
%	\Biggr.
%	\\
%	\Biggl.
	-
	\left(
		\beta_{n}^{-1} \E[X_{n1}]
		+
		R_{n}
	\right)^{2}
\Biggr\}
(1 + o(1))
\end{multline*}
where
\begin{eqnarray*}
R_{n} = 
 \begin{cases}
		C_{2}\eps_{n} \beta_{n}^{-\min(1,\alpha)} \tau^{1-\alpha}, &\text{if $\alpha\ne 1$,}\\
		C_{3}\eps_{n} \beta_{n}^{-1} \log(\beta_{n}), &\text{if $\alpha=1$,}
\end{cases}
\end{eqnarray*}
and \(C_{1}, C_{2}, C_{3}>0.\) % Note that under our choice of \(\beta_{n}\), \(n \eps_{n} \beta_{n}^{-\alpha} \to c_{1}\). 
We conclude that if \(n \beta_{n}^{-2} \to 0\) (that is, \(\gamma_{2}<1-\alpha/2\)), then \(\lim_{\tau\to 0} \lim_{n \to \infty} [ G_{n} ] =0;\) otherwise the last limit is infinite. 
At the same time, \eqref{an} yields for \(\alpha \ne 1,\)% for any \(H>0\)
\begin{eqnarray}\label{ann1}
a_{n}&=& n \left[
	\frac{
		1-\eps_{n}
	}
	{\beta_{n}}
	\E \left[
		X_{n1}
	\right]
+
\frac{
\eps_{n}
}
{\beta_{n}^{\min(1,\alpha)}
}
%H^{1-\alpha} 
\frac{\alpha}{|1-\alpha|}
\right] + o(1).
\end{eqnarray}
For instance, if \(\alpha<1, \) then 
\begin{eqnarray*}
a_{n}=
\frac{
n \E \left[
X_{n1}
\right]
}
{
 \beta_{n}
}+
\frac{\alpha}{1-\alpha}
+o(1),
\end{eqnarray*}
where the first summand in the r.h.s. is of the order \(n^{1-(1-\gamma_{2})/\alpha}\). Therefore, the choice of  \(a_{n}\) differs in the cases \(\gamma_{2} \in (1-\alpha, 1-\alpha/2)\), \(\gamma_{2}=1-\alpha\), and \(\gamma_{2} \in(0, 1- \alpha)\), and this observation leads to different choices of \(\kappa_{n} = a_{n}\beta_{n},\) see \eqref{kappan}. %In the case \(\alpha>1\), we get from \eqref{ann1}, \(a_{n} = n \E[X_{n1}]/\beta_{n}.\) 
Finally, in the case \(\alpha=1,\)
\begin{eqnarray*}
a_{n}&=& n \left[
	\frac{
		1-\eps_{n}
	}
	{\beta_{n}}
	\E \left[
		X_{n1}
	\right]
+
\frac{
\eps_{n} \alpha
}
{\beta_{n}
}
\log(\beta_{n})
\right]
%H^{1-\alpha} 
+ o(1) = \frac{n \E[Z_{n1}]}{\beta_{n}}+o(1),
\end{eqnarray*}
where we use \eqref{EZ}.
\item \underline{\(\gamma_{2} \in \left(
	1-\alpha/2, \; (2-\alpha) \gamma_{1}
\right).\)} In this case,  we take \(\beta_{n} = \sqrt{n \Var(X_{n1})}\). Under this choice, 
the conditions (a) and (b) from Part 3 of Proposition~\ref{con} hold. 
The choice \(a_{n}=n \E[X_{n1}]/\beta_{n}\) follows from \eqref{ann1}.
 \item \underline{\(\gamma_{2} < 
	1-\alpha \gamma_{1}.\)} It is easy to see that the infinite smallness condition \eqref{inf} is not fulfilled. Note that this case was considered separately in Theorem~\ref{thm3} (i). 

\end{enumerate}
\textbf{Proof of Theorem~\ref{thm1} (ii) and (iii).}  The proof directly follows from the application of the well-known Slutsky theorem. For instance, Theorem~\ref{thm3}(ii) yields that 
\begin{eqnarray*}
\frac{S_{n} - n \E\left[ X_{n1}\right]
	}
	{
n \E\left[ X_{n1}\right]
	}
=\frac{S_{n} - n \E\left[ X_{n1}\right]
	}
	{
		\sqrt{n \Var(X_{n1})}
	}
	\cdot
	\frac{\sqrt{\Var(X_{n1})}}{\sqrt{n}\E\left[ X_{n1}\right]} \eqp 0,
\end{eqnarray*}
since the first  multiplier tends in distribution to the standard normal law, and the second tends to 0.

\section{Acknowledgment} The author is grateful to Prof. Stanislav Molchanov (UNC Charlotte, USA, and Higher School of Economics, Moscow, Russia)  for the supervision of  this research.

\bibliographystyle{plain}
\bibliography{Panov_bibliography}

\end{document}